\documentclass[a4paper]{article}

\usepackage{amsmath,amssymb,amsthm}

\newcommand{\dd}{\displaystyle}
\newcommand{\surg}{\Sigma_g}
\newcommand{\sur}{\Sigma}
\newcommand{\surgd}{\Sigma_g^{\ud}}
\newcommand{\surdq}{\Sigma_g^{(\ud)}}

\def\OO{{\mathcal{O}}}
\def\CC{{\mathbb C}}
\def\RR{{\mathbb R}}
\def\ZZ{{\mathbb Z}}
\def\NN{{\mathbb N}}
\def\TCC{{\underline{\mathbb C}}}
\def\DD{{\mathfrak D}}
\def\JJ{{\mathcal{J}}}
\def\SS{{\mathfrak{S}}}

\def\MM{{\mathcal{M}}}
\def\MP{{\mathcal{P}}}
\def\TT{{\mathcal{T}}}

\def\arf{{\mathcal{A}}}

\def\DB{{\bar{\partial}}}

\def\CCC{{\mathcal{C}}}

\DeclareMathOperator{\ddet}{Det}
\DeclareMathOperator{\rk}{rk}
\DeclareMathOperator{\rank}{rk}
\DeclareMathOperator{\rest}{rest}
\DeclareMathOperator{\e}{e}
\DeclareMathOperator{\coker}{coker}
\DeclareMathOperator{\PD}{PD}

\newcommand{\riso}[1][N]{\RR Aut(#1)}
\newcommand{\raut}[1][N]{\RR GL(#1)}
\newcommand{\rsaut}[1][N]{\RR SL (#1)}

\newcommand{\Det}[1][N]{\ddet(#1)}
\newcommand{\rop}[1][N]{\RR \CCC (#1)}
\newcommand{\ropdx}[1][N]{\RR \CCC_{D_{\ux}} (#1)}
\newcommand{\ropc}[1][N]{\RR \CCC_{compat} (#1)}
\newcommand{\ropcp}[1][N]{\RR \CCC_{compat}^+ (#1)}

\newcommand{\ppm}{\mathfrak{p}^{\pm}}
\newcommand{\pp}{\mathfrak{p}^{+}}
\newcommand{\mmp}{\mathfrak{p}^{-}}
\newcommand{\PP}[1][N]{\mathfrak{p}^+(#1)}
\newcommand{\PPX}{\mathfrak{p}^+_X}
\newcommand{\dethun}{\det(H^1(\surg,\RR)_{-1})}
\newcommand{\dethzero}{\RR_{w}}
\newcommand{\dethunww}{H^1_w(\RR\surg,\RR)}

\newcommand{\ud}{\underline{\mathbf{d}}}
\newcommand{\ux}{\underline{\mathbf{x}}}
\newcommand{\uz}{\underline{\mathbf{z}}}

\newtheorem{Definition}{Definition}
\newtheorem{Lemma}{Lemma}
\newtheorem{Theorem}{Theorem}
\newtheorem{Proposition}{Proposition}
\newtheorem{Corollary}{Corollary}
\newtheorem{Remark}{Remark}

\title{Real bundle automorphisms, Cauchy Riemann operators and orientability of moduli spaces}
\author{R\'emi Cr\'etois\footnote{Research is supported in part by the ERC project TROPGEO} \\ \small{Universit\'e de Gen\`eve} \\ \small{remi.cretois@unige.ch}}

\begin{document}

\maketitle

\centerline{\textbf{Abstract}}
This paper consists in a very brief English summary of the results appearing in French in \cite{article1} and \cite{article2}. We will omit the proofs and focus on explaining our approach and theorems. This paper is not intended to be published. Questions are welcome.

In \cite{article1} and \cite{article2} we consider a complex vector bundle $N$ equipped with a real structure $c_N$ over a real curve $(\surg,c_\sur)$ of genus $g$. We compute the sign of the action of an automorphism of $(N,c_N)$ on the orientations of the determinant line bundle over the space of Cauchy-Riemann operators on $(N,c_N)$. 

We first consider the automorphisms lifting the identity on $\surg$ (in \cite{article1}). In this case, we obtain the sign as a product of two terms. The first one computes the signature of the permutations induced by the automorphisms acting in the $Pin^\pm$ structures on the real part of $(N,c_N)$. The second one comes from the action of the automorphisms of $(N,c_N)$ on the bordism classes of real $Spin$ structures on $(\surg,c_\sur)$. We then study the general case in \cite{article2}. As an application of these results, we compute the first Stiefel-Whitney class of the moduli space of real pseudoholomorphic curves in many cases.
\\

\noindent
\textsc{Classification AMS 2010}: 14H60, 53D45

\section{Setting}

Let $\surg$ be a connected, oriented and closed surface of genus $g$ and $c_\sur$ a real structure on $\surg$, i.e. an orientation-reversing involutive diffeomorphism of $\surg$. We call such a pair $(\surg,c_\sur)$ a real curve. Fix a complex vector bundle $\pi : N\rightarrow\surg$ and equip it with a real structure $c_N$, i.e. $c_N$ is an involutive anti-$\CC$-linear automorphism of $N$ which lifts $c_\sur$. The fixed point set $\RR\surg$ of $c_\sur$ is called the real part of the curve. It is a $1$-dimensional submanifold of $\surg$, which can be empty, and over which the fixed point set $\RR N$ of $c_N$ forms a real vector bundle. When $c_\sur$ or $c_N$ induce an involution on a vector space or on a module, we will denote with a subscript $+1$ (resp. $-1$) the eigenspace associated to the $+1$ (resp. $-1$) eigenvalue of this involution.

Real curves are classified up to diffeomorphism by two invariants: the number $0\leq k\leq g+1$ of connected components of $\RR \surg$, and whether $\surg\setminus \RR\surg$ is connected or not. In the latter case, we say that the curve is a separating curve. Given a real curve, the complex vector bundles equipped with real structures are classified up to isomorphism by their rank, first Chern class and by the first Stiefel-Whitney class of their real part.

Fix integers $1 \ll q \ll l$ and a real number $p > 1$ such that $qp > 2$. Denote by $\RR J (\surg)$ the set of all complex structures $J$ of class $\CCC^l$ on $\surg$ which are compatible with the orientation and such that $c_\sur$ is $J$-antiholomorphic. It is a contractible Banach manifold (voir \cite{wel1}).

\begin{Definition}
A real Cauchy-Riemann operator on $(N,c_N)$ is a pair $(\DB, J)$ consisting of a complex structure $J\in\RR J(\surg)$ and a $\CC$-linear operator
\[
\DB : L^{q,p}(\surg,N)\rightarrow L^{q-1,p}(\surg,\Lambda^{0,1}\surg\otimes N),
\]
which is $c_N$-equivariant and which satisfies Leibniz's rule:
\[
\forall f \in \CCC^{\infty}(\surg,\CC),\ \forall v\in L^{k,p}(\surg,N),\ \DB(fv) = f\DB(v) + \DB_J(f)\otimes v,
\]
where $\DB_J = \frac{1}{2}(d + i\circ d\circ J)$.
\end{Definition}

A real Cauchy-Riemann operator $\DB$ on $(N,c_N)$ induces a Fredholm operator from $L^{q,p}(\surg,N)_{+1}$ to $L^{q-1,p}(\surg,\Lambda^{0,1}\surg\otimes N)_{+1}$. We will write $H^0_{\DB}(\surg,N)_{+1}$ and $H_{\DB}^1(\surg,N)_{+1}$ for its kernel and cokernel.

The set $\rop$ of all real Cauchy-Riemann operators on $(N,c_N)$ is a contractible space. One can define a real line bundle $\Det$ over it whose fibre over an operator $\DB$ is its determinant line $\ddet(\DB) = \Lambda^{\max}_{\RR}H^0_{\DB}(\surg,N)_{+1}\otimes \left(\Lambda^{\max}H^1_{\DB}(\surg,N)_{+1}\right)^*$ (see \cite{MDS}, Theorem A.2.2). This bundle is orientable but not canonically oriented. To understand this lack of canonical orientation, we study the action of the automorphism group of $(N,c_N)$ on the orientations of $\Det$.

\begin{Definition}
An automorphism of $(N,c_N)$ is a pair $(\Phi,\varphi)$ consisting of
\begin{enumerate}
\item an orientation-preserving diffeomorphism $\varphi : (\surg,c_\sur)\rightarrow (\surg,c_\sur)$ which is $\ZZ/2\ZZ$-equivariant,
\item a $\ZZ/2\ZZ$-equivariant map $\Phi : (N,c_N) \rightarrow (N,c_N)$ which is $\CC$-linear in the fibers, and such that $\pi\circ\Phi = \varphi\circ\phi$.
\end{enumerate}
We let $\riso$ be the group of all automorphisms of $(N,c_N)$.
\end{Definition}

The group $\riso$ acts on the orientations of $\Det$ in the following way. First, it acts on sections of $N$ and on $(0,1)$-forms with value in $N$ on one hand and on the elements of $\RR J(\surg)$ on the other hand like this:
\[
(\Phi,\varphi)_* s = \Phi(s)\circ\varphi^{-1},
(\Phi,\varphi)^*\alpha = \Phi^{-1}\circ\alpha\circ d\varphi,
(\Phi,\varphi)^*J = d\varphi^{-1}\circ J\circ d \varphi,
\]
where $s$ is a section of $N$, $\alpha$ is a $(0,1)$-form with value in $N$ and $J\in\RR J(\surg)$. This induces an action on $\rop$
\[
(\Phi,\varphi)^*(\DB,J) = ((\Phi,\varphi)^*(\DB((\Phi,\varphi)_*)),(\Phi,\varphi)^*J),
\]
and then on the bundle $\Det$.

Our first goal was to compute the action of an element of $\riso$ on the orientations of $\Det$ in topological terms. We state the results we obtained in this direction (\cite{article1} and \cite{article2}) in \S \ref{resultats1}.

The question raised above was motivated by the orientability problem of the moduli spaces of real pseudoholomorphic curves in real symplectic manifolds. More precisely, consider a real symplectic manifold $(X,\omega,c_X)$, i.e. a symplectic manifold $(X,\omega)$ equipped with an anti-symplectic involution $c_X$. We denote by $\RR J_{\omega}(X)$ the set of all almost complex structures $J$ on $X$ which are $\omega$-compatible and such that $c_X$ is $J$-antiholomorphic. Given an integer $g$ and a class $d\in H_2(X,\ZZ)$, we let $\pi: \RR \MM_g^d(X)\rightarrow \RR J_{\omega}(X)$ by the universal moduli space of real pseudoholomorphic curves of genus $g$ realizing the class $d$ and which are somewhere injective. We consider only somewhere injective curves to avoid some analytical problems caused by multiply covered curves. Indeed, under this hypothese, the fiber $\RR \MM_g^d(X,J) = \pi^{-1}\{J\}$ is a smooth finite dimensional manifold for all generic $J\in\RR J_{\omega}(X)$. Moreover, given an element $u\in\RR \MM_g^d(X,J)$ which is an immersed curve, there is a monodromy morphism $\mu : \pi_1(\RR\MM_g^d(X,J))\rightarrow \pi_0(\riso[N_u])$, where $N_u$ is the normal bundle to $u$. Then, for each $\gamma \in \pi_1(\RR \MM_g^d(X,J))$, the first Stiefel-Whitney class of $\RR\MM_g^d(X,J)$ computed against $\gamma$ is zero if and only if $\mu(\gamma)$ acts trivially on the orientations of $\Det[N_u]$.

We state the results we obtained in \cite{article2} concerning this second problem in \S \ref{resultats2}. Let us note that similar results have been obtained recently by Georgieva and Zinger in \cite{georgieva}, but using different methods.

\section{Real automorphisms of a bundle and the determinant bundle}\label{resultats1}

Let $(\surg,c_\sur)$ be a real curve and $(N,c_N)$ be a complex vector bundle equipped with a real structure over this curve. A first step in our study is to describe the stucture of the automorphism group $\riso$ of $(N,c_N)$. To that end, we introduce two subgroups of this group:
\begin{itemize}
\item the group $\raut$ of automorphisms which lift the identity of $\surg$, 
\item the subgroup $\rsaut\subset\raut$ of elements $f\in \raut$ such that the determinant automorphism $\det(f) : \det(N)\rightarrow \det(N)$ is the identity.
\end{itemize}

Here, $\det(N)$ with no capital letter is the standard determinant of a vector bundle : $\det(N) = \Lambda_{\CC}^{\rk N} N$.

The following group of diffeomorphisms of $\surg$ will also appear in the following Lemma: for $w\in H^1(\RR\surg,\ZZ/2\ZZ)$, we write
\[
\RR Diff^+(\surg, w) = \{\phi : \surg\rightarrow \surg\ \text{ orientation-preseving diffeomorphism } |\ \phi\circ c_\sur = c_\sur \circ \phi,\ \phi^*w = w\}.
\]

\begin{Lemma}[Lemma 2.2 of \cite{article1}]\label{seq}
The following sequences are exact:
\begin{equation}
  \label{suite1}
  0\rightarrow \raut \rightarrow \riso \rightarrow \RR Diff^+(\surg,w_1(\RR N))\rightarrow 0
\end{equation}
\begin{equation}\label{suite2}
    0 \rightarrow \rsaut \longrightarrow \raut \xrightarrow{\det} \RR\CCC^{\infty}(\surg,\CC^*) \rightarrow 0
  \end{equation}
  \begin{equation}\label{suite3}
    0 \rightarrow \ker(\rest) \longrightarrow \rsaut \xrightarrow{\rest} \pi_0(SL(\RR N)) \rightarrow 0
  \end{equation}
  where:
  \[
  \begin{array}{c c c c}
    \rest : &\rsaut & \rightarrow & \pi_0(SL(\RR N)) \\
    &\Phi      & \rightarrow & [\Phi_{|\RR N}].
  \end{array}
  \]
Moreover, $\ker(\rest)$ is connected.
\end{Lemma}

As the action of an element of $\riso$ on the orientations of $\Det$ depends only on its homotopy class, the Lemma \ref{seq} suggest that we proceed in three steps: first we consider only the subgroup $\rsaut$, then we consider $\raut$ and finally we go on to the general case. As we want to express this action in topological terms, we introduce the necessary topological objects we need before stating the different results.

\subsection{$Pin$ structures: the case of $\rsaut$}

Suppose first that $\RR\surg\neq \emptyset$ and let us recall that we put $k = b_0(\RR \surg)$. An element $\Phi$ of $\riso$ acts by pullback on the $(\ZZ/2\ZZ)^k$-principal space $Pin^{\pm}(\RR N)$ consisting of the $Pin^{\pm}$ structures on each connected component of $\RR N$. Thus, we get permutations $\Phi_{\ppm}$ of $Pin^{\pm}(\RR N)$. One must be careful as the two permutations $\Phi_{\pp}$ and $\Phi_{\mmp}$ are not the same in general. More precisely, given a diffeomorphism $\varphi\in\RR Diff^+(\surg,w_1(\RR N))$ it induces a permutation $\sigma_{\varphi}^-$ of the $2k_-$ elements set consisting of the orientations of the connected components of $\RR\surg$ over which $\RR N$ is not orientable. Then, denoting by $\varepsilon(\sigma)$ the signature of any permutation $\sigma$ of a finite set, we get the following.

\begin{Proposition}[Proposition 1.1 of \cite{article2}]
  Let $(N,c_N)$ be a complex vector bundle equipped with a real structure over a real curve $(\surg,c_\sur)$ with non empty real locus. Let $(\Phi,\varphi)\in\riso$. Then $\varepsilon(\Phi_{\pp}) = \varepsilon(\Phi_{\mmp})\varepsilon(\sigma_{\varphi}^-)$.
\end{Proposition}

If we restrict our attention to the elements of $\rsaut$, then their action on $Pin^{\pm}(\RR N)$ is enough to describe their homotopy class modulo $2$. This will be enough for our purpose because the action of an element of $\rsaut$ on the orientations of $\Det$ depends only on its class in $\pi_0(\rsaut)/2\pi_0(\rsaut)$. This last remark is also valid when $\RR\surg = \emptyset$ in the sense that in that case $\rsaut$ is connected, so all its elements act trivially on the orientations of $\Det$.

We can now state our first result.

\begin{Theorem}[Theorem 3.1 of \cite{article1}]\label{thp}
 Let $(N,c_N)$ be a complex vector bundle equipped with a real structure over a real curve $(\surg,c_\sur)$. The action of an element $\Phi\in\rsaut$ on the orientations of $\Det$ is given by $\varepsilon(\Phi_{\ppm})$.
\end{Theorem}

\subsection{Real $Spin$ structures: the case of $\raut$}

From the exact sequence (\ref{suite2}) in Lemma \ref{seq}, we have to understand the structure of the group $\raut[\det(N)]$ which is canonically equal to $\RR\CCC^{\infty}(\surg,\CC^*)$. In fact, the Theorem \ref{thp} allows us to restrict our attention to the case where $N$ is of rank one. Indeed, one has the following: 

\begin{Proposition}[Lemma 4.1 of \cite{article1}]\label{ranko}
  Let $(N,c_N)$ be a complex vector bundle equipped with a real structure over a real curve $(\surg,c_\sur)$. The action of an element $\Phi\in\raut$ on the orientations of $\Det$ is equal to the product of $\varepsilon(\Phi_{\ppm})$ by the action of $\det(\Phi)$ on the orientations of $\Det[\det(N)]$.
\end{Proposition}

Let us suppose that $\rank(N) = 1$, $\RR \surg\neq \emptyset$, and fix a complex structure $J\in\RR J(\surg)$. We then consider the set $\RR Spin(\surg)$ of all real $Spin$ structures on $(\surg,c_\sur)$, that is the set of square roots of $K_\sur$ admitting a real structure. The cardinal of $\RR Spin(\surg)$ is $2^{g+k-1}$.

This set can be partitioned using, for every $w\in H^1(\RR\surg,\ZZ/2\ZZ)$ such that $w([\RR\surg]) = g+1 \mod 2$, the set $\RR Spin(\surg,w)$ of square roots of $K_\sur$ admitting a real structure whose real locus is a real vector bundle having $w$ as first Stiefel-Whitney class. We will then say that an element of $\RR Spin(\surg,w)$ is a real $Spin$ structure with first Stiefel-Whitney class $w$.

\begin{Definition}
We say that two real $Spin$ structures are in the same bordism class if they have the same first Stiefel-Whitney class and if they are in the same bordism class in the usual sense (see \cite{kirby}).
\end{Definition}

The group $\raut$ acts on the set of real $Spin$ structures via the canonical correspondance between the elements of $\raut$ and of $\raut[T\surg]$. Moreover, it can be proved that $\raut$ acts on the bordism classes of real $Spin$ structures of a given first Stiefel-Whitney class (see Lemma 4.9 of \cite{article1}).

\begin{Definition}
  For all $w\in H^1(\RR\surg,\ZZ/2\ZZ)$ such that $w([\RR\surg]) = g+1 \mod 2$, we let $\arf^w : \raut \rightarrow \ZZ/2\ZZ$ be the morphism corresponding to the action of $\raut$ on the bordism classes of real $Spin$ structures whose first Stiefel-Whitney class is $w$.
\end{Definition}

This action is important in our results, however, it is not enough to detect completely the action of $\raut$ on the orientations of $\Det$. For example, the automorphism $-1$ acts trivially on $Spin$ structures, but its action on the orientations of $\Det$ is given by the Riemann-Roch theorem, i.e. it acts trivially if and only if $\deg(N) + 1 - g = 0 \mod 2$. In fact, we will see later that the morphism $\arf^{w_1(\RR N)}$ coincides with the action of $\raut$ on the orientations of $\Det$ only when $\deg(N) + 1 - g = 0\mod 2$.

To treat the general case, for each $w \in H^1(\surg,\ZZ/2\ZZ)$ we let $\beta^w : \raut\rightarrow \ZZ/2\ZZ$ be the morphism computing the product of the signs of the determinant of an element of $\raut$ over each component of $\RR\surg$.

\begin{Definition}[See Lemma 4.11 of \cite{article1}]
For each $f\in\raut$, let
  \begin{itemize}
  \item  $s_{top}(f) = \beta^w(f)\arf^w(f) \in\ZZ/2\ZZ$, for any $w\in H^1(\surg,\ZZ/2\ZZ)$ such that $w([\RR\surg]) = g+1\mod 2$
  \item  $s_N(f) = s_{top}(f)\beta^{w_1(\RR N)}(f)\in\ZZ/2\ZZ$.
  \end{itemize}
\end{Definition}

We can extend these definitions to the case where $\rank(N)>1$ by first taking the determinant of the automorphism $f\in\raut$ then applying the previous definitions.

\begin{Remark}
  \begin{itemize}
  \item In the case where $\deg(N) = g+1 \mod 2$, then we can take $w = w_1(\RR N)$ in the definition of $s_{top}(f)$, so we get $s_N(f) = \arf^{w_1(\RR N)}(f)$.  \item In the case where the curve is separating, i.e. when $\surg\setminus \RR\surg$ has two connected components, one can check that the sign $s_{top}(f)$ is always $1$, so we get $s_N(f) = \beta^{w_1(\RR N)}(f)$. In this case, the sign $s_N(f)$ boils down to the product of the actions of $f$ on the orientations of the orientable components of $\RR N$.
  \end{itemize}
\end{Remark}

We can now state the following.

\begin{Theorem}[Theorem 4.3 of \cite{article1}]\label{thsp}
 Let $(N,c_N)$ be a complex vector bundle equipped with a real structure over a real curve $(\surg,c_\sur)$ with non-empty real locus. The action of an element $\Phi\in\raut$ on the orientations of $\Det$ is given by $\varepsilon(\Phi_{\ppm})s_N(\Phi)$.
\end{Theorem}

To prove this theorem, we use a well-chosen family of automorphisms and compute the action of each of them. In fact, we have the following Proposition.

\begin{Proposition}[Proposition 1.8 of \cite{article2}]
  Let $(N,c_N)$ be a complex line bundle equipped with a real structure over a real curve $(\surg,c_\sur)$. Let $a\subset \surg$ be a closed oriented simple curve. Consider the following automorphisms $f_a\in\raut$:
  \begin{enumerate}
    \item If $a$ is a connected component of $\RR\surg$, then choose a tubular neighborhood of $a$ of the form $(\theta,t)\in S^1\times [-1,1]$, where $a$ corresponds to $t=0$ and $c_\sur(\theta,t) = (\theta,-t)$. Set $f_a(\theta,t) = -\e^{i\pi t}$ and extend it by $1$ outside of this neighborood.
    \item If $a$ is globally fixed by $c_\sur$, then choose a tubular neighborhood of $a$ of the form $(\theta,t)\in S^1\times [-1,1]$, where $a$ corresponds to $t=0$ and $c_\sur(\theta,t) = (-\theta,-t)$. Set $f_a(\theta,t) = -\e^{i\pi t}$ and extend it by $1$ outside of this neighborood.
    \item If $a\cap c_\sur(a) = \emptyset$, then choose a tubular neighborhood of $a$ disjoint from its conjugate and of the form $(\theta,t)\in S^1\times [-1,1]$, where $a$ corresponds to $t=0$. Set $f_a(\theta,t) = -\e^{i\pi t}$ on this neighborhood and extend it by $\overline{f_a\circ c_\sur}$ on its conjugate and by $1$ everywhere else.
  \end{enumerate}
In the first case, $f_a$ preserves the orientations of $\Det$ if and only if $\RR N$ is non-orientable over $a$. In the second case, $f_a$ never preserves the orientations of $\Det$. In the last case, $f_a$ always preserves the orientations of $\Det$.
\end{Proposition}

Note the in the previous Proposition, we do not assume the real curve to have a non-empty real locus.

As an application of the Theorem \ref{thsp}, we study the determinant bundle over the real Picard group and compute its first Stiefel-Whitney class (see Theorem 4.4 of \cite{article1}). We obtain the same results as Okonek and Teleman in \cite{oktel}.

\subsection{Using Divisors: the general case}

The Proposition \ref{ranko} extends nicely to the general case. Indeed, consider the trivial rank $n$ bundle $(\TCC^{\oplus n},conj)$ with its standard real structure over $(\surg,c_\sur)$. For each $\varphi\in\RR Diff^+(\surg)$, we define $\Phi_{\varphi}\in\riso[\TCC^{\oplus n}]$ by
\[
\begin{array}{c c c c}
  \Phi_\varphi : & \TCC^{\oplus n} & \rightarrow & \TCC^{\oplus n} \\
  & (x,v) & \mapsto & (\varphi(x),v).
\end{array}
\]

\begin{Proposition}[Proposition 1.2 of \cite{article2}]\label{rankoo}
 Let $(N,c_N)$ be a complex vector bundle equipped with a real structure over a real curve $(\surg,c_\sur)$. The action of an element $(\Phi,\varphi)\in\riso$ on the orientations of $\Det$ is given by the product of $\varepsilon(\Phi_{\pp})$ with the action of $\det(\Phi)$ on the orientations of $\Det[\det(N)]$ and with the action of $\Phi_{\varphi}\in\riso[\TCC^{\rank(N)-1}]$ on the orientations of $\Det[\TCC^{\rank(N)-1}]$.
\end{Proposition}

Compared to Proposition \ref{ranko}, there is one additional term which we discuss now.

\begin{Proposition}[Proposition 1.3 of \cite{article2}]\label{trivial}
  Let $n$ be a positive integer and $\varphi \in \RR Diff^+(\surg)$. The action of $\Phi_{\varphi}\in\riso[\TCC^{\oplus n}]$ on the orientations of $\Det[\TCC^{\oplus n}]$ is given by $\det(\varphi^*)^n$ where $\varphi^* : H^1(\surg,\RR)_{-1}\rightarrow H^1(\surg,\RR)_{-1}$ is induced by $\varphi$.
\end{Proposition}

There remains only to study the action of $\det(\Phi)$ on the orientations of $\Det[\det(N)]$ to conclude the general case. That is why we suppose in the following that $N$ is of rank one.

Unfortunately, the action on real $Spin$ structures introduced in the previous paragraph does not extend well to the whole group $\riso$. Indeed, on one hand, in general there is no natural way to split the exact sequence (\ref{suite1}) in Lemma \ref{seq}. On the other hand, to be able to extend the above action, we need a way to separate the behaviour of an element of $\riso$ in the fiber of $N$ from its behaviour on $\surg$. More precisely, we need to understand how a element of $\riso$ ``twists'' the fibers of $N$, i.e. some analog of the index map $\raut \rightarrow Hom(H_1(\surg,\ZZ),H_1(\CC^*,\ZZ))$.

To this end, we introduce some new objects.

\begin{Definition}
  A $c_\sur$-invariant divisor $D = \dd\sum_{i} a_i x_i$, with $a_i\in\ZZ$ and $x_i\in\surg$ is said to be compatible with $(N,c_N)$ if and only if
  \begin{itemize}
  \item it is of degree $\deg(N)$,
  \item its restriction to $\RR\surg$ is Poincar\'e dual to $w_1(\RR N)$.
  \end{itemize}
\end{Definition}

Given a divisor $D$ compatible with $(N,c_N)$ and an element $(\DB,J)\in\rop$ such that $(N,\DB)\cong \OO_{\surg,J}(D)$, the set of real meromorphic sections of $(N,\DB)$ with divisor $D$ is a real line. Now we can understand how an element $(\Phi,\varphi)$ of $\riso$ such that $\varphi_* D = D$ behaves in the fibers of $N$ by comparing a real meromorphic section of $(N,\DB)$ with divisor $D$ to its image by $\Phi$. It is with this idea in mind that we will proceed in the following.

Let $\ud$ be a $4$-tuple of positive integers $(r^+,r^-,s^+,s^-)\in \NN^4$. Let $\surgd$ be the subset of $\left(\RR\surg\right)^{r^+ + r^-}\times \left(\surg\right)^{s^++s^-}\setminus \Delta$, where $\Delta$ is the diagonal, consisting of elements $\ux = (\ux^+,\ux^-,\uz^+, \uz^-)$, where $\uz^{\pm} =(z_1^{\pm},c_\sur(z_1^{\pm}),\ldots,z_{s^{\pm}}^{\pm}, c_\sur(z_{s^{\pm}}^{\pm}))$ such that the divisor $D_{\ux} = \dd\sum_{x\in\ux^+} x - \dd\sum_{x\in\ux^-} x + \dd\sum_{i=1}^{s^+} (z_i^+ + c_\sur(z_i^+)) - \dd\sum_{i=1}^{s^-} (z_i^- + c_\sur(z_i^-))$ is compatible with $(N,c_N)$.

\begin{Definition}
  If $\surgd$ is non-empty, we say that $\ud$ is adapted to $(N,c_N)$.
\end{Definition}

Consider the bundle $\pi : (pr_1)^*\rop \rightarrow \RR J(\surg)\times\surgd$, where $pr_1: \RR J(\surg)\times\surgd\rightarrow \RR J(\surg)$ is the first projection. Define the subset $\ropc\subset (pr_1)^*\rop$ whose fiber over $\ux\in\surgd$ is
\[
\ropdx = \{(\DB,J)\in\rop\ |\ (N,\DB)\cong \OO_{\surg,J}(D_{\ux})\}.
\]

We will say that an element $\DB$ of $\ropdx$ is polarised if we have chosen an orientation of the real line made of the real meromorphic sections of $(N,\DB)$ having $D_{\ux}$ as divisor.

The set $\ropc$ is a $\raut/\RR^*$-principal bundle over $\RR J(\surg)\times\surgd$, and its double cover $\ropcp$ consisting of polarised operators is a $\raut/\RR^*_+$-principal bundle over the same space. Unfortunately, the group $\riso$ acts on those two spaces in a way which is not compatible with the bundle structure. That is why we use the following definition.

\begin{Definition}
  Let $\raut^+\subset\raut$ be the subgroup consisting of automorphisms acting trivially on the orientations of $\Det$. Fix $\ud\in\NN^4$ adapted to $(N,c_N)$.
  \begin{itemize}
  \item If the curve $(\surg,c_\sur)$ is separating and all the components of $\RR N$ are non-orientable, let $\DD_{\ud}(N)$ be the trivial $\ZZ/2\ZZ$-principal bundle over $\RR J(\surg)\times\surgd$.
  \item In all other cases, let $\DD_{\ud}(N)$ be the $\ZZ/2\ZZ$-principal bundle $\ropcp/\raut^+$.
  \end{itemize}
\end{Definition}

Let us write $\SS_{2,s^{\pm}}$ the group generated by the elements of the permutation group $\SS_{2s^{\pm}}$ of the form $(2i - 1\ 2i)$, for $1\leq i\leq s^{\pm}$, and those of the form $2i-1\mapsto 2\sigma(i)-1$ and $2i\mapsto 2\sigma(i)$, for $\sigma\in\SS_{s^{\pm}}$. The group $\SS_{\ud} = \SS_{r^+}\times \SS_{r^-}\times \SS_{2,s^+}\times \SS_{2,s^-}$ acts freely on $\surgd$ by renumbering the points. Let $\surdq$ be the quotient of $\surgd$ by $\SS_{\ud}$. The action of $\SS_{\ud}$ lifts naturally to $\DD_{\ud}(N)$. The quotient $\DD_{(\ud)}(N)$ is a $\ZZ/2\ZZ$-principal bundle over $\RR J(\surg)\times\surgd$.

One can show that the bundle $\DD_{\ud}(N)$ is trivialisable, but that in general it is not the case for $\DD_{(\ud)}$ (see Proposition 1.9 of \cite{article2}). More precisely, if we consider the real vector bundle $\RR\JJ_{(\ud)}$ over $\RR J(\surg)\times\surdq$ whose fiber over $\{J\}\times (\ux)$ is the vector space generated by the elements of $(\ux)$, then the bundle $\DD_{(\ud)}(N)\otimes \det(\RR\JJ_{(\ud)})$ is orientable. In particular, when the curve is non separating, the bundle $\DD_{(\ud)}(N)$ will be non-orientable over a loop in $\surdq$ consisting in exchanging two complex conjugated points. 

Given an element $(\Phi,\varphi)\in\riso$, this bundle $\DD_{\ud}(N)$ will allow us two separate the contribution to the action of $\Phi$ on $\Det$ coming from $\varphi$ and that coming from the fiber part of $\Phi$.

Of course, when the rank of $N$ is greater than $1$, we extend the previous definitions by taking the determinant of $N$.

\begin{Theorem}[Theorem 1.2 of \cite{article2}]\label{bigres}
  Let $(N,c_N)$ be a complex vector bundle equipped with a real structure over a real curve $(\surg,c_\sur)$, and let $\ud$ be a $4$-tuple adapted to $(N,c_N)$. The determinant bundle $Det$ over $\rop\times\surdq$ is canonically isomomorphic to
\[
\PP\otimes \DD_{(\ud)}(N)\otimes \RR \JJ_{(\ud)}\otimes T_{(\ud)} \otimes \det(H^1(\surg,\RR)_{-1})^{\otimes \rank(N)},
\]
where $\PP$ is the trivial bundle with fiber $\dd\bigotimes_{i=1}^k Pin^+((\RR N)_i)$ and $T_{(\ud)}$ is the line bundle whose fiber over $(\DB, (\ux)) \in\rop\times\surdq$ is $\dd\bigotimes_{x\in (\ux^+)}T^*_x\RR\surg$.
\end{Theorem}

What we mean here by canonical is that we construct an isomorphism which commutes with the action of $\riso$ on both sides.

The proof of the Theorem \ref{bigres} uses real elementary negative transformations on $N$, which allows us to replace $N$ with another bundle of a different degree, namely the trivial bundle, while keeping track of the automorphisms and their action on the determinant bundle.

The result given by the Theorem \ref{bigres} is not very easy to use because of the bundle $\DD_{(\ud)}(N)$ which remains a bit mysterious. That is why we also treat two important particular cases.

\subsubsection{The separating case}

Let $(N,c_N)$ be a complex vector bundle equipped with a real structure over a real curve $(\surg,c_\sur)$. Suppose that the curve is separating, i.e. that $\surg\setminus\RR \surg$ is non-connected. Denote by $\RR\surg^-$ the union of connected components of $\RR\surg$ over which $\RR N$ is non-orientable and set $k_- = b_0(\RR \surg^-)$. Let $O_{\RR N}$ be the real line $\dd\bigotimes_{c\subset \RR\surg\setminus\RR\surg^-} o((\RR N)_{|c})$ where $o((\RR N)_{|c})$ is the real line generated by the two orientations of $\RR N$ over $c$. Lastly, let $H$ be the real line generated by the two complex orientations of $\RR \surg$.

\begin{Corollary}[Corollary 1.4 of \cite{article2}]\label{sepa}
Let $(N,c_N)$ be a complex vector bundle equipped with a real structure over a real curve $(\surg,c_\sur)$, and suppose that the curve is separating. Then there is a canonical isomorphism between $\Det$ and the real line bundle
\[
\PP\otimes O_{\RR N}\otimes \det(H^0(\RR\surg^-,\RR))\otimes \det(H^1(\surg,\RR)_{-1})^{\otimes \rank{N}} \otimes H^{\otimes\frac{\deg(N)+k_-}{2}}.
\]
\end{Corollary}

\subsubsection{The $Spin$ case}

Let $(N,c_N)$ be a complex vector bundle equipped with a real structure over a real curve $(\surg,c_\sur)$ with non-empty real locus. Suppose that $N$ is of even degree and that $\RR N$ is orientable. Then we consider the set $\RR Spin(N)$ of real $Spin$ structures on $(N,c_N)$. Fixing a real Cauchy-Riemann operator on $(N,c_N)$, a real $Spin$ structure on $(N,c_N)$ is equivalent to a holomorphic square root $L$ of $\det(N)$ admitting a real structure. There are two choices of a real structure on $L$ and each one induces an orientation on $\RR N$ opposite of the other. Moreover, those two orientation do not depend on the initial choice of Cauchy-Riemann operator but only on the $Spin$ structure. Thus, to each element $\xi\in\RR Spin(N)$ we associate the line $O^{\xi}_{\RR N}$ generated by the two orientations of $\RR N$ given by $\xi$. We can also define the first Stiefel-Whitney class $w_{\xi}\in H^1(\RR\surg,\ZZ/2\ZZ)$ of $\xi\in\RR Spin(N)$ as $w_1(\RR L)$ for any of the two real structure on $L$. Lastly, for all $w\in H^1(\RR\surg,\ZZ/2\ZZ)$, we set $H_w^1(\RR\surg,\RR) = \dd\bigotimes_{w([\RR\surg]_i) = 1}H^1((\RR\surg)_i,\RR)$.

\begin{Corollary}[Corollary 1.5 of \cite{article2}]\label{spin}
  Let $(N,c_N)$ be a complex vector bundle equipped with a real structure over a real curve $(\surg,c_\sur)$ with non-empty real locus. Suppose that $N$ is of even degree and that $\RR N$ is orientable. Then there is a canonical isomorphism
\[
\Det = \PP\otimes \det(H^1(\surg,\RR)_{-1})^{\otimes \rank(N)}\otimes (O^{\xi}_{\RR N})^{\otimes 1-g}\otimes H_{w_{\xi}}^1(\RR\surg,\RR)
\]
of line bundles over $\rop\times\RR Spin(N)$.
\end{Corollary}

\section{Orientability of Moduli Spaces}\label{resultats2}

\subsection{Real Teichm\"uller space}

Suppose in this section that the curve $(\surg,c_\sur)$ is of genus at least $2$. We give a result concerning the action of a diffeomorphism of $(\surg,c_\sur)$ on the orientations of the real Teichm\"uller space of the curve.

\begin{Definition}
The real Teichm\"uller space associated to the real curve $(\surg,c_\sur)$ is the quotient of $\RR J(\surg)$ by the action of the subgroup $\RR Diff_0(\surg)$ of $\RR Diff^+(\surg)$ consisting of those diffeomorphisms which are homotopic to the identity,
\[
\TT(\surg,c_\sur)= \RR J(\surg)/ \RR Diff_0(\surg).
\]
\end{Definition}

As a consequence of the Proposition \ref{trivial} we have the following Theorem.

\begin{Theorem}[Theorem 1.1 of \cite{article2}]
  Let $(\surg,c_\sur)$ be a real curve of genus at least $2$. The orientation bundle of $\TT(\surg,c_\sur)$ is canonically isomorphic to the trivial bundle $\det(H^1(\surg,\RR)_{-1})$.
\end{Theorem}

In particular, one sees that the moduli spaces obtained by quotienting the Teichm\"uller spaces by the diffeotopy group will not always be orientable.

\subsection{Moduli spaces of real curves in real symplectic manifolds}

From now on, let $(X,\omega,c_X)$ be a real symplectic manifold of dimension $2n$ at least $4$. We will give the results on the orientability of moduli spaces we obtained using the Theorem \ref{bigres} and the Corollaries \ref{sepa} and \ref{spin}. But first we recall some necessary definitions.

Fix integers $1 \ll q \ll l$ and a real $p > 2$. Recall that $J_{\omega}(X)$ is the space of all almost complex structures on $X$ of class $C^l$ which are tamed by $\omega$, and $\RR J_{\omega}(X)$ is the subset of those for which $c_X$ is anti-holomorphic. Both are Banach manifolds and contractible.

\begin{Definition}
  A parametrized pseudo-holomorphic curve of genus $g$ in $X$ is a triple $(u,J_\sur,J)\in L^{q,p}(\surg,X)\times J(\surg)\times J_{\omega}(X)$ such that
\[
d u + J\circ d u \circ J_\sur = 0.
\]
A pseudo-holomorphic curve is in the class $d\in H_2(X,\ZZ)$ if $u_*([\surg]) = d$. It is said to be somewhere injective if there exists a non empty open subset $U\subset \surg$ such that $u_{|U}$ is an immersion and
\[
\forall x\in U, u^{-1}\{u(x)\} = \{x\}.
\]
For $d\in H_2(X,\ZZ)$, we set
\[
\begin{array}{l l }
  \MP^d_{g,r}(X) = \{(u,J_\sur,J,\uz)\in L^{q,p}(\surg,X)\times J(\surg)\times J_{\omega}(X)\times \surg^r\ | & d u + J\circ d u \circ J_\sur = 0 \\
  & u_*([\surg]) = d \\
& u\text{ is somewhere injective}\\
& \forall 1\leq i\neq j\leq r,\ z_i\neq z_j\}.
\end{array}
\]
\end{Definition}

Fix a permutation $\tau$ of $\{1,\ldots,r\}$ of order $2$. Using $\tau$, the group $Diff(\surg)$ of diffeomorphisms of $\surg$ of class $C^{l+1}$ acts on $\MP^d_{g,r}(X)$ by reparameterization:
\[
\begin{split}
(\varphi,(u,J_\sur,J,\underline{z}))\in Diff(\surg)\times\MP_{g,r}^d(X)\hspace{7cm} \\ \mapsto \left\{\begin{array}{l r}\left(u\circ\varphi^{-1},(\varphi^{-1})^*J_\sur,J,\varphi(\underline{z})\right) & \text{if } \varphi\in Diff^+(\surg)\\
\left(c_X\circ u\circ\varphi^{-1},-(\varphi^{-1})^*J_\sur,J,(\varphi(z_{\tau(1)}),\ldots,\varphi(z_{\tau(r)}))\right) & \text{if not}.
\end{array}\right.
\end{split}
\]

If an element of $\MP_{g,r}^d(X)$ is fixed by a diffeomorphism, then this diffeomorphism has to be a real structure $c_\sur$ on $\surg$ and is unique. We then write $\RR_{c_\sur}\MP_{g,r}^d(X)$ the set of curves which are fixed by $c_\sur$, and $\RR_{\tau}\MP^d_{g,r}(X)$ all the fixed points of $\MP^d_{g,r}(X)$ under the action of $Diff(\surg)$. Those two sets are Banach manifolds of class $C^{l-q}$.

We set $\RR_{\tau}\MM_{g,r}^d(X)$ to be the quotient of $\MP_{g,r}^d(X)$ by the free action of $Diff(\surg)$. It is again a Banach manifold of class $C^{l-q}$. Moreover, there is a forgetful map $\Pi : \RR_{\tau}\MM^d_{g,r}(X) \rightarrow \RR \MM^d_{g}(X)$ and a Fredholm map $\pi : \RR\MM_{g}^d(X)\rightarrow \RR J_{\omega} (X)$. We are interested in the orientability of the fibers of this last map. More precisely, there is a real line bundle $\ddet(\pi)$ over $\RR \MM_g^d(X)$ whose fiber over $[u,J_\sur,J]$ is $\Lambda^{\max}\ker(d_{[u,J_{\sur},J]} \pi) \otimes \left(\Lambda^{\max}\coker(d_{[u,J_{\sur},J]} \pi)\right)^*$, and its restriction to each regular fiber is canonically isomorphic to the orientation bundle of the fiber. It is in some sense the bundle of relative orientations of $\pi : \RR \MM^d_g(X)\rightarrow \RR J_{\omega}(X)$.

To take into account the contribution of marked points, we introduce the tautological line bundle $L_r$ over $\RR_\tau \MM_{g,r}^d(X)$ whose fiber over $[u,J_\sur,J,\uz]$ with real structure $c_\sur$ is $\det(\RR_{c_\sur}\dd\bigoplus_{z\in\uz} T_z\surg)$.

Lastly we will use the bundles $\PPX$ and $\dethun$ over $\RR \MM^d_g(X)$ whose fibers over $[u,J_\sur,J,\uz]$ with real structure $c_\sur$ are respectively $\dd\bigotimes_{i=1}^{b_0(\RR\surg)}Pin^+((\RR E_u)_i)$, where $E_u = u^*TX$, and $\dethun$.

In the next three paragraphs, we give the results we obtained in \cite{article2} concerning this bundle.

\subsubsection{Separating curves}

We begin with a particular case which was already studied (see \cite{wel4}, \cite{puignau}, \cite{fooo} and \cite{Solomon}). We partition $\RR_\tau\MM^d_{g,r}(X)$ into two components: $\RR_\tau^{sep}\MM_{g,r}^d(X)$ containing the separating curves and $\RR_\tau^{nsep}\MM_{g,r}^d(X)$ containing the non-separating curves.

First, let us note that there is a locally constant function $k_- : \RR_\tau\MM^d_{g,r}(X)\rightarrow \{0,\ldots,g+1\}$ which counts for every $[u,J_\sur,J,\uz]$ the number of connected components of the real part of $E_u = u^*TX$ that are non-orientable. Then, we consider the real line bundle $O_X$ over $\RR_\tau\MM_{g,r}^d(X)$ whose fiber over $[u,J_\sur,J,\uz]$ is $O_{\RR E_u}$.

Lastly, for a separating curve $[u,J_\sur,J,\uz]$ with real structure $c_\sur$, we set $H_{c_\sur}$ to be the real line generated by the two complex orientations of $\RR\surg$ and $\RR\surg^-$ be the reunion of the components of $\RR\surg$ over which $\RR E_u$ is non-orientable. The line bundle $H$ and $\dethzero$ over $\RR_\tau^{sep}\MM_{g,r}^d(X)$ are then the line bundles with respective fiber $H_{c_\sur}$ and $\dethzero$ over $[u,J_\sur,J,\uz]$.

The following Theorem follows essentially from the Corollary \ref{sepa}.

\begin{Theorem}[Corollary 2.1 of \cite{article2}]\label{thsep}
  The line bundle $\Pi^*\ddet(\pi)$ over $\RR^{sep}_\tau\MM_{g,r}^d(X)$ is canonically isomorphic to
\[
  \dethun^{\otimes n-1}\otimes \dethzero\otimes
  H^{\otimes \frac{c_1(X)d+k_-}{2}} \otimes \PPX \otimes O_X.
\]
\end{Theorem}

In this case, we can also show that the line bundles $L_r$ and $H$ are orientable when $\tau$ is not the identity, i.e. when there is at least one pair of complex conjugated marked points.

\begin{Corollary}[Corollary 2.2 of \cite{article2}]
  Let $(X,\omega,c_X)$ be a symplectic manifold of dimension $2n$ at least $4$. Assume that $\RR X$ is $Spin$. Then for all $d\in H_2(X,\ZZ)$, for all $r\geq 3$, for all permutation $\tau$ of $\{1,\ldots,r\}$ of order $2$ with at least one fixed point, and for all $J\in\RR J_{\omega}(X)$ generic enough, the manifold $\RR_{\tau}\MM_{0,r}^d(X,J) = \pi^{-1}\{J\}$ is orientable.

More generally, assume that $n$ is odd. Then for all $d\in H_2(X,\ZZ)$, for all $g\geq 0$, $r\geq 2$, for all non-trivial permutation $\tau$ of $\{1,\ldots,r\}$ of order $2$, and for all $J\in\RR J_{\omega}(X)$ generic enough, the manifold $\RR_{\tau}^{sep}\MM_{g,r}^d(X,J) = \pi^{-1}\{J\}$ is orientable.
\end{Corollary}

One can also derive from Theorem \ref{thsep} other more particular results: here we used the $Spin$ hypothese to kill the contributions from $\dethzero$, $\PPX$ and $O_X$, and the hypothese on the permutations to kill the contributions of $H$ and of the marked points.

\subsubsection{$Spin$ case}

Take $d\in H_2(X,\ZZ)$ such that $c_1(X)d$ is even and let $\RR_{\tau}\MM^{d,0}_{g,r}(X)$ be the reunion of the components of $\RR_{\tau}\MM^{d}_{g,r}(X)$ containing the curves $[u,J_\sur,J,\uz]$ that have non-empty real locus and such that $w_1(\RR E_u) = 0$. Consider $\RR_{\tau}\MM^{d,Spin}_{g,r}(X)$ the cover of degree $2^{g+k-1}$ over $\RR_{\tau}\MM^{d,0}_{g,r}(X)$ associated to the $(\ZZ/2\ZZ)^{g+k-1}$-principal bundle whose fiber over $[u,J_\sur,J,\uz]$ consists of the real $Spin$ structures on $E_u$. On the one hand there is a real line bundle $O_X^{Spin}$ over $\RR_{\tau}\MM^{d,Spin}_{g,r}(X)$ associated to the $\ZZ/2\ZZ$-principal bundle whose fiber over $([u,J_\sur,J,\uz],\xi)$ consists of the two semi-orientations of $\RR E_u$ coming from $\xi$. On the other hand, all the real $Spin$ structure in a given connected component of $\RR_{\tau}\MM^{d,Spin}_{g,r}(X)$ have the same first Stiefel-Whitney class, so we can define the bundle $\dethunww$ over $\RR_{\tau}\MM^{d,Spin}_{g,r}(X)$ with fiber $\dd\bigotimes_{w([\RR\surg]_i) = 1}H^1((\RR\surg)_i,\RR)$ over $([u,J_\sur,J,\uz],\xi)$, where $\xi$ has $w$ as first Stiefel-Whitney class.

The Theorem \ref{thspin} follows essentially from the Corollary \ref{spin}.

\begin{Theorem}[Corollary 2.3 of \cite{article2}]\label{thspin}
Let $d\in H_2(X,\ZZ)$ such that $c_1(X)d$ is even. The bundle $\Pi^*\ddet(\pi)$ over $\RR_\tau\MM_{g,r}^{d,Spin}(X)$ is canonically isomorphic to
\[
\dethun^{\otimes n-1}\otimes \dethunww \otimes(O_X^{Spin})^{\otimes 1-g}\otimes \PPX.
\]
\end{Theorem}

\begin{Corollary}[Corollary 2.4 of \cite{article2}]
  Let $X_\delta$ be a smooth hypersurface of degree $\delta$ in $\CC P^N$, $N\geq 4$, $N= 0$ or $3\mod 4$. Assume that $\RR X_{\delta}$ is non-empty and that $\delta = N+1\mod 4$ and $\delta \leq N+1$. Let $d\in H_2(X_\delta,\ZZ)$, $r\geq 1$ and $\tau$ be a permutation of $\{1,\ldots, r\}$ of order $2$ with at least one fixed point. Then for all $J\in\RR J_{\omega}(X)$ generic enough,
\[
w_1(\RR_{\tau}\MM^d_{g,r}(X_\delta,J)) = w_1(L_r) + (\delta-1)w_1(\dethun).
\]
\end{Corollary}

Here, we used the hypotheses on the degree to guarantee that $X_\delta$ has a real $Spin$ structure whose first Stiefel-Whitney class is trivial, and that its real part is $Spin$. This way we get a section of $\RR_{\tau}\MM^{d,Spin}_{g,r}(X_\delta,J)\rightarrow\RR_{\tau}\MM^d_{g,r}(X_\delta,J)$, and we can kill $\PPX$, $O_X^{Spin}$ and $\dethunww$.

\subsubsection{Polarisations}

We now go on to the general case. To this end, we will consider polarisations on $(X,\omega,c_X)$.

\begin{Definition}
  A real divisor on $(X,\omega,c_X)$ is a finite sum $D = \sum a_V V$, with $a_V\in\ZZ$ and $V\subset X$ is a smooth symplectic submanifold of dimension $2n-2$, such that $(c_X)_*D = D$.
\end{Definition}

Given such a real divisor $D$ on $(X,\omega,c_X)$, one can consider the subset $\RR J_D(X)$ of $\RR J_{\omega}(X)$ consisting of the almost complex structures $J\in\RR J_{\omega}(X)$ such that all the components of $D$ are $J$-holomorphic. When there exists such a structure which is moreover compatible with $\omega$, then one can show that $\RR J_D(X)$ is a contractible Banach manifold.

\begin{Definition}
  A real divisor $D$ is said to be admissible if
  \begin{itemize}
  \item all its components intersect one another transversally,
  \item all its components have multiplicity $\pm 1$,
  \item there exists an element in $\RR J_D(X)$ which is compatible with $\omega$.
  \end{itemize}
We then denote by $D_r$ the part of $D$ containing the components which are stable under $c_X$.

A polarisation of $(X,\omega,c_X)$ is an admissible divisor $D$ such that $D^{\PD} = c_1(X)\in H^2(X,\ZZ)$ and $(\RR D_r)^{\PD} = w_1(\RR X)\in H^1(\RR X,\ZZ/2\ZZ)$.

A polarizing section associated to $D$ is a real section $s$ of $\det_J(TX) = \Lambda^n_J TX$ for some $J\in\RR J_D(X)$ such that
\begin{itemize}
\item $s^{-1}(\{0\}) = D$,
\item $s$ vanishes transversally along $D$ except where two components intersect,
\item the sign of each component of $D$ is given by comparing the orientation of the component coming from $s$ and that coming from $\omega$.
\end{itemize}
\end{Definition}

Given a polarisation $D$ on $(X,\omega,c_X)$, we set
\[
\begin{array}{l l }
\RR\MM_g^d(X,D)^{\pitchfork} = \{[u,J_\sur,J]\in\RR\MM_g^d(X)\ | & J\in \RR J_D(X)  \\
& u(\surg)\not\subset D\\
& \forall V\subset D,\ u(\surg)\pitchfork V \\
&\forall V\neq V'\subset D,\ u(\surg)\cap V \cap V' = \emptyset\}.
\end{array}
\]
In other words, we look at curves which intersect $D$ transversally and along a single component at a time. The space $\RR\MM_g^d(X,D)^{\pitchfork}$ is an open subset of the space of $J$-holomorphic curves, for $J\in \RR J_D(X)$, not contained in $D$, which is a Banach manifold. Moreover, writing $D = \sum a_V V$, there is a line bundle $T_D$ defined over $\RR\MM_g^d(X,D)^{\pitchfork}$ whose fiber over $[u,J_\sur,J]$ is given by $\RR\JJ_{u^{*}(D)}\otimes \dd\bigotimes_{
  \begin{subarray}{c}
    x \in u^{-1}(V)\\
    a_V = 1
  \end{subarray}} T_x^*\RR\surg$.

\begin{Theorem}[Theorem 2.2 of \cite{article2}]\label{thgen}
  Let $(X,\omega,c_X)$ be a real symplectic manifold polarised by the divisor $D$ given by a polarizing section. Then the bundle $\ddet(\pi: \RR\MM_g^d(X,D)^{\pitchfork}\rightarrow \RR J_D(X))$ is canonically isomorphic to
\[
\PPX\otimes \dethun^{\otimes n-1}\otimes T_D.
\]
\end{Theorem}

The Theorem \ref{thgen} is the translation of the Theorem \ref{bigres} in the context of moduli spaces. Here, the bundle $\DD_{(\ud)}(N)$ does not appear thanks to the polarizing section.

In \cite{article2} we also give some results about the existence of polarisations. In particular, we show that any real symplectic manifold admits a polarisation given by a polarizing section. We also give some concrete examples in the case of smooth projective hypersurfaces.

\bibliographystyle{plain}
\bibliography{article}

\end{document}